\documentclass[11pt]{article}
\usepackage{amsmath}
\usepackage{graphicx}
\usepackage{epsfig}
\usepackage{amsfonts}
\usepackage{amssymb}
\usepackage{placeins}
\usepackage{cases}
\usepackage[margin=1in]{geometry}
\usepackage{authblk}
\usepackage[titletoc,toc,title]{appendix}
\usepackage{epstopdf}
\usepackage{changes}
\usepackage{xcolor}
\usepackage{bm}
\usepackage{enumerate}

\usepackage{accents}

\title{Data-Driven Model Identification Near a Supercritical Hopf Bifurcation Using Phase-Based Approaches}

\begin{document}
\author[1]{Dan Wilson \thanks{corresponding author:~dwilso81@utk.edu}}
\affil[1]{Department of Electrical Engineering and Computer Science, University of Tennessee, Knoxville, TN 37996, USA}
\maketitle

\begin{abstract}
A data-driven model identification strategy is developed for dynamical systems near a supercritical Hopf bifurcation with nonautonomous inputs.  This strategy draws on phase-amplitude reduction techniques, leveraging an analytical representation for the phase and amplitude response curves of the Hopf normal form to infer system parameters.  Fitting can be performed by recording the system output during the relaxation to the stable limit cycle after applying as few as two carefully timed pulse inputs.  This strategy is illustrated in two examples with relevance to circadian oscillations.  In each example, the proposed model identification strategy allows for the formulation, solution, and implementation of a closed loop nonlinear optimal control problem.
\end{abstract}

\section{Introduction} \label{intsec}

Phase-based reduction techniques have a rich history in the analysis and control of weakly perturbed oscillators \cite{kura84}, \cite{erme10}, \cite{piet19}, \cite{wils22review}.  This approach allows for the transformation of a general system of equations
\begin{equation}  \label{dynsintro}
\dot{x} = F(x,u),
\end{equation}
where $x \in \mathbb{R}^N$ is the state, $u \in \mathbb{R}^M$ is an input, and $F$ gives the dynamics to a phase reduced form
\begin{equation} \label{predintro}
    \dot{\theta} = \omega + Z(\theta) u,
\end{equation}
where $\theta \in \mathbb{S}$, $\omega$ is the unperturbed natural frequency, and $Z$ is an appropriately sized vector that captures the response to inputs as a function of phase.  The transformation of \eqref{dynsintro} to \eqref{predintro} enables the rigorous mathematical analysis of the dynamics of both coupled and externally forced oscillations in applications that would otherwise be intractable \cite{erme10}, \cite{wils22review}.   The relatively simple form of \eqref{predintro} allows for a straightforward data-driven model identification strategy:~by applying a brief pulse of input and considering the resulting change in phase one can obtain a pointwise measurement of $Z(\theta)$.  Repeating this procedure multiple times by applying inputs at different phases, one can obtain a data-driven estimate of $Z(\theta)$.  This so-called `direct method' \cite{gala05}, \cite{neto12} has been used in a variety of experimental applications to study behaviors such as synchronization, phase locking, and entrainment \cite{gala05}, \cite{khal03}, \cite{nabi132}, \cite{wils15optent}, \cite{snar15}.  

While phase reduction and phase models of the form \eqref{predintro} are remarkably useful, a fundamental limitation is that it is only valid in the limit that inputs are infinitesimally small, i.e.,~its accuracy is only guaranteed in a weakly perturbed regime.  Practically, the allowable magnitude of the input is related to the nonunity Floquet multipliers that govern the decay in directions transverse to the periodic orbit.  Recent years have seen a flurry of interest in the development of phase-based reduction methods that can accommodate stronger magnitude inputs \cite{lets20}, \cite{cast20}, \cite{wils20highacc}, \cite{geng20}, \cite{wils21adapt}, \cite{kure13}, \cite{kure22}.  While such techniques provide a more accurate representation of the dynamics in applications where strong magnitude inputs must be applied, it is typically more difficult to infer the associated reduced order models from data.  Associated model identification methods have been developed for some of these approaches \cite{wils23direct}, \cite{wils20ddred}, \cite{wils18operat}, but they typically require much more data (and subsequently many more experiments) than what is required for the fitting of a standard phase model of the form \eqref{predintro}.

In this work, a direct method approach for model identification is developed and investigated for systems near a supercritical Hopf bifurcation.   This strategy draws on phase-amplitude reduction techniques, first obtaining an analytical representation for the phase and amplitude response curves of the Hopf normal form and using this information to uniquely fit system parameters.  The resulting model is valid in a close vicinity of the Hopf bifurcation, rather than just in a close vicinity of the unperturbed limit cycle like the standard phase reduction of the form \eqref{predintro}.  The proposed approach can be implemented by recording a general output during the relaxation to the limit cycle after applying as few as two carefully timed pulse perturbations.  The organization of this paper is as follows:~Section \ref{backsec} gives necessary background information about phase-amplitude reduction techniques using isochrons and Floquet coordinates and discusses associated methods for model identification using the direct method.  Section \ref{hopfsec} considers the normal form of the Hopf bifurcation in relation to these phase-amplitude reduction techniques.  A careful consideration of the resulting reduced order equations yields a model identification strategy that can infer all terms of the Hopf normal form with a combination of passive observations and a small number of pulse perturbations.  Section \ref{exampsec} provides illustrative examples on two models with relevance to circadian oscillations.  In both examples, the resulting model identification strategy enables the application of closed-loop optimal control algorithms.  Section \ref{concsec} provides concluding remarks.

\section{Background} \label{backsec}

\subsection{Phase Reduction}

Consider an ordinary differential equation of the form 
\begin{align} \label{maineq}
\dot{x} &= F(x,u), \nonumber \\
    y &= g(x),
\end{align}
where $x \in \mathbb{R}^N$ is the state, $y \in \mathbb{R}^K$ is the output, $u \in \mathbb{R}^M$ is an input, $F$ gives the dynamics, and $g$ maps the state to the output.  Suppose that when $u$ is held constant at $u = 0$, Equation \eqref{maineq} has a stable $T$-periodic orbit $x^\gamma$.  Rather than considering the state of \eqref{maineq}, it can be useful to consider the timing of oscillations.  In this case, one can define a phase $\theta \in [0,2\pi]$ for $x \in x^\gamma$.  By convention, $\theta$ is typically scaled so that $\frac{d \theta}{dt} = \frac{2 \pi}{T} = \omega$ when taking $u = 0$.  Isochrons can be used to define phase in the basin of attraction of the limit cycle \cite{guck75}, \cite{winf01}.  Subsequently changing to phase coordinates, one can represent Equation \eqref{maineq} in terms of its phase dynamics \cite{erme10}
\begin{equation} \label{pred}
    \dot{\theta} = \omega + Z(\theta)u,
\end{equation}
where $Z(\theta) = \frac{\partial \theta}{\partial x}^T \frac{\partial F}{\partial u}$ with all partial derivatives evaluated at $x^\gamma(\theta)$, and $^T$ denotes the transpose.  Equation \eqref{pred} is valid in the limit that $u$ is small.  Numerically, $Z(\theta)$ can be computed by first finding $\partial \theta/\partial x$, which has solutions that satisfy the adjoint equation \cite{brow04}, \cite{erme91}
\begin{equation} \label{adjeq}
    \frac{d \nabla \theta}{dt} = -J^T \nabla \theta.
\end{equation}
Above, $\nabla \theta \equiv \frac{\partial \theta}{\partial x}$ and $J$ denotes the Jacobian of $F$, both evaluated on the periodic orbit at $\theta(t)$.  In addition to satisfying \eqref{adjeq}, $\nabla \theta$ must also be scaled so that $\nabla \theta ^T F(x^\gamma) = \omega$.

\subsection{Phase-Amplitude Reduction}
Phase reduction of the form \eqref{pred} only considers the timing of oscillations.  To incorporate amplitude-based effects, one can augment the phase equation with Floquet coordinates.  Towards this end, first consider the linear approximation of \eqref{maineq} near the periodic orbit 
\begin{equation} \label{lineq}
    \dot{x} = J \Delta x,
\end{equation}  
where $\Delta x = x - x_{p}^\gamma(\theta)$ and $J$ is the Jacobian evaluated at $x^\gamma(\theta(t))$.  Noting that $J$ is $T$-periodic and leveraging Floquet theory \cite{jord07}, provided the monodromy matrix of this linear time varying system is diagonalizable, Equation \eqref{lineq} admits solutions of the following form
\begin{equation} \label{floqform}
x - x^\gamma(\theta) = \sum_{j = 1}^{N-1} \psi_j g_j(\theta).
\end{equation}
Above, $g_j(\theta,p) \in \mathbb{C}^N$ is a Floquet eigenfunction associated with the Floquet coordinate $\psi_j \in \mathbb{C}$.  Above, the contribution from the Floquet eigenfunction with associated Floquet multiplier $\lambda_N = 1$ has been absorbed by the phase coordinate to give $N-1$ total Floquet coordinates.  Note that each $g_j$ is unique to a constant scaling.   Considering the underlying equation \eqref{maineq} and changing to Floquet coordinates, the dynamics of the Floquet coordinates are \cite{wils16isos}
\begin{align} \label{isoeq}
\dot{\psi}_j &= \kappa_j \psi_j + I_j(\theta) u, \nonumber \\
j &= 1,\dots,N-1,
\end{align}
where $I_j(\theta,p) = \frac{\partial \psi_j}{\partial x}^T  \frac{\partial F}{\partial u}$ with all partial derivatives evaluated at $x_{p}^\gamma(\theta)$, and $\kappa_j \in \mathbb{C}$ is the Floquet exponent associated with the $j^{\rm th}$ Floquet eigenfunction. Equations \eqref{pred} and \eqref{isoeq} together constitute a phase-amplitude-based representation for the underlying system \eqref{maineq} that is valid provided $\Delta x$ remains small.   In many cases, it can be useful to truncate the most rapidly decaying $\psi_j$ (as gauged by the magnitude of the associated $\kappa_j$) thereby arriving at a phase amplitude reduction.  Indeed, the inclusion of amplitude coordinates allows for the formulation and solution of control problems that phase reduction alone cannot accommodate \cite{wils16isos}, \cite{mong19b}, \cite{taka21}.  Numerically, each $I_j(\theta)$ can be computed by first finding $\frac{\partial \psi_j}{\partial \theta}$ which has solutions that satisfy \cite{wils16isos}
\begin{equation} \label{isoeq2}
    \frac{d \nabla \psi_j}{dt} = -(J^T - \kappa_j {\rm Id}) \nabla \psi_j,
\end{equation}
where $\nabla \psi_j \equiv \frac{\partial \psi_j}{\partial \theta}$ evaluated on the periodic orbit at $\theta(t)$ and ${\rm Id}$ is an appropriately sized identity matrix. Additionally, the solution to \eqref{isoeq2} must be scaled so that $\nabla \psi_j^T g_j(\theta) = 1$ for all states on the periodic orbit.  

Note that the phase and amplitude reduced equations from \eqref{pred} and \eqref{isoeq} are only valid in the limit that the input magnitude is small relative to the magnitude of the nonunity Floquet multipliers.  Floquet coordinates can be generalized to the entire basin of attraction of the limit cycle using the concept of isostable coordinates, which are level stets of principle Koopman eigenfunctions \cite{mezi20}, in order to improve the accuracy of the phase amplitude reduction when larger magnitude inputs are required \cite{wils20highacc}.  For a broader discussion regarding the use of isostable coordinates in control applications, the interested reader is referred to \cite{wils22review}.

\subsection{Direct Methods for Inference of Phase and Amplitude Response Curves} \label{directback}
In an experimental setting, the equations that comprise $F$ are often unavailable thereby precluding the use of Equations \eqref{adjeq} and \eqref{isoeq2} for finding $Z(\theta)$ and $I_j(\theta)$, respectively.  Alternatively, one can infer $Z(\theta)$ and from data by applying pulse inputs and observing the subsequent relaxation to the limit cycle. To summarize approaches from \cite{neto12} and \cite{gala05}, for simplicity of exposition, consider $M = 1$ so that $u \in \mathbb{R}$ (although this is not a strict requirement).  First, it is necessary to define a Poincar\'e section $\Gamma_0$ constructed as a $K-1$ dimensional hyperplane in the space of outputs that is transversal to $g(x^\gamma(\theta))$ at $\theta = 0$.   When $g(x) \in \Gamma_0$, $\theta$ is approximately zero with the mismatch depending on the local curvature of the isochrons relative to the $\Gamma_0$ hyperplane.  Next, the frequency and associated period ($\omega$ and $T$, respectively) of the system \eqref{maineq} must be approximated from passive observations.

Starting with estimation of $Z(\theta)$, for an initial condition $x = x^\gamma(\theta_0)$ one can apply a short pulse of input $u(t) = m$ lasting $t_0$ time units.  By determining the resulting change in phase, for instance, as gauged by the subsequent crossing of $\Gamma_0$ relative to the expected crossing had the pulse input not been applied, one can obtain an estimate of the phase response curve according to 
\begin{equation} \label{phasedirect}
Z(\theta_0) \approx \frac{\Delta \theta}{m t_0}.
\end{equation}
This procedure can be repeated at different choices of $\theta_0$ allowing $Z(\theta)$ to subsequently be fit to a Fourier series basis $Z(\theta) = b_0 + \sum_{k = 1}^\nu \left[ a_k \sin(\theta) + b_k \cos(\theta) \right]$ choosing $\nu$ appropriately to avoid overfitting.  

For the amplitude coordinate dynamics, suppose that $|\kappa_j| = O(1/\epsilon)$ for $j \geq 2$ where $0 < \epsilon \ll 1$, i.e.,~so that all but one Floquet coordinate decays rapidly.  In this case, it is possible to use the same pulses inputs to obtain an estimate for the slowly decaying Floquet exponent, $\kappa_1$, and the corresponding response curve, $I_1(\theta)$.  A careful description of this strategy is given in \cite{wils18operat} and summarized here.  Let $\theta(t)$ at some time $t = t_p$ be equal to $\theta_0$.  Applying a pulse input $u(t) = m$ starting at $t = t_p$ and lasting $t_0$ time units, one can define $t^*_k$ to be the $k^{\rm th}$ crossing of the $\Gamma_0$ Poincar\'e section after the application of the pulse input.  Subsequently letting $\tau_k = t^*_k - t^*_{k-1}$ the slope of a linear regression of $k$ versus $\log(\tau_k-T)/T$ is equal to $\kappa_1$.  For the same pulse input, as shown in \cite{wils18operat}, the amplitude response curve can be computed according to
\begin{equation} \label{floqdirect}
    I_1(\theta_0) = \frac{\rho}{m t_0}  \frac{2 \pi (k-1) - \omega(t^*_k - t^*_1)}{  \exp(\kappa_1(t_k^* - t_p)) - \exp(\kappa_1(t_1^*-t_p))},
\end{equation}
for  and $k \geq 2$.  Above, $\rho$ is a constant that depends on the underlying system \eqref{maineq}.  As with the phase response curve, this procedure can be repeated for different choices of $\theta_0$ allowing $I_1(\theta)$ to subsequently be fit to a Fourier basis.  Note that data from the same pulse input can be used to obtain a pointwise estimate of both $Z(\theta_0)$ and $I_1(\theta_0)$.  Related strategies for estimating $Z(\theta)$ and $I_k(\theta)$ for systems with multiple non-negligible Floquet coordinates have been developed leveraging machine learning techniques \cite{ahme24}.  References \cite{mong19}, \cite{wils18operat}, and \cite{wils23direct} provide a more complete description of these model identification strategies for the interested reader.  In general these approaches require the collection of enough data to fully resolve both $Z(\theta)$ and $I_1(\theta)$, potentially requiring a substantial of trials.

\section{Data-Driven Model Identification for Hopf Oscillators} \label{hopfsec}

As shown here, for oscillatory systems of the form \eqref{maineq} that result from a supercritical Hopf bifurcation, the phase-amplitude coordinate framework enables model identification by observing the relaxation to the limit cycle following as few as two pulse perturbations.  The following analysis will assume that $M = 1$ so that $u \in \mathbb{R}$, but appropriate modifications could be made straightforwardly to consider the case where $M>1$.  

\subsection{Normal Form for a Hopf Bifurcation With Applied Control}
The autonomous form of the Hopf bifurcation has been widely studied \cite{wigg03}, \cite{guck83}.  This basic structure is adapted here for use with an additional nonautonomous forcing term.   To begin suppose that $F$ from \eqref{maineq} is at least $C^5$ differentiable with an unstable fixed point $x_0$ with $F(x_0,0) = 0$. Suppose also that the system is close to a supercritical Hopf bifurcation so that the linearized vector field has a complex-conjugate pair of eigenvalues $\lambda_{1,2}$ with ${\rm Real}(\lambda_{1,2}) \approx 0$ with the remaining eigenvalues $\lambda_j$ having $|{\rm Real}(\lambda_j)|$ sufficiently bounded away from 0 for $j \geq 3$.    Note that $u$ is not explicitly assumed to be a bifurcation parameter.  Considering the normal form of the Hopf bifurcation, when taking $u = 0$, it is possible to transform \eqref{maineq} to the following form:
\begin{align} \label{hnf}
    \dot{\tilde{x}} &= \tilde{\alpha} \tilde{x} - \tilde{\beta} \tilde{y} + (\tilde{a} \tilde{x} - \tilde{b} \tilde{y})(\tilde{x}^2 + \tilde{y}^2), \nonumber \\
    \dot{\tilde{y}} &= \tilde{\beta} \tilde{x} + \tilde{\alpha} \tilde{y} + (\tilde{b} \tilde{x} + \tilde{a} \tilde{y})(\tilde{x}^2 + \tilde{y}^2),
\end{align}
which is valid up to fifth order of accuracy.  Above, $\tilde{x} \in \mathbb{R}$ and $\tilde{y} \in \mathbb{R}$ represent Cartesian coordinates on the center manifold associated with the Hopf bifurcation that result from the analytic coordinate change $\tilde{x} = h_1(x)$, $\tilde{y} = h_2(x)$ and $\tilde{\alpha}, \tilde{\beta}, \tilde{a},\tilde{b} \in \mathbb{R}$ are coefficients associated with the Hopf normal form.   In radial coordinates, i.e.,~for which $\tilde{x} = \tilde{r} \sin(\tilde{\varphi})$ and $\tilde{y} = \tilde{r} \cos(\tilde{\varphi})$ the dynamics are $\dot{\tilde{r}} = \alpha \tilde{r} + a \tilde{r}^3  + O(\tilde{r}^5)$, $\dot{\tilde{\varphi}} = \beta + br^2 + O(\tilde{r}^4)$.  Focusing on the dynamics of \eqref{hnf} $u = 0$, one can write
\begin{align} \label{hopfcontrolpre}
    \dot{\tilde{x}} &= \frac{\partial h_1}{\partial x} \cdot \frac{d x}{dt} \nonumber \\
    &  =  \frac{\partial h_1}{\partial x} \cdot \bigg( F(x,0) + \frac{\partial F}{\partial u} u + O(u^2)\bigg) \nonumber \\
    &  = \tilde{\alpha} \tilde{x} - \tilde{\beta} \tilde{y} + (\tilde{a}\tilde{x} - \tilde{b} \tilde{y})(\tilde{x}^2 + \tilde{y}^2) +  H_1 u + O(u^2) + O(u |x - x_0|), \nonumber \\
    \dot{\tilde{y}} &= \frac{\partial h_2}{\partial x} \cdot \frac{d x}{dt} \nonumber \\
    &  =  \frac{\partial h_2}{\partial x} \cdot \bigg( F(x,0) + \frac{\partial F}{\partial u} u + O(u^2) \bigg) \nonumber \\
    &  = \tilde{\beta} \tilde{x} + \tilde{\alpha} \tilde{y} + (\tilde{b}\tilde{x} + \tilde{a} \tilde{y})(\tilde{x}^2 + \tilde{y}^2) + H_2 u + O(u^2) + O(u |x - x_0|),
\end{align}
where the `dot' denotes the dot product and $H_i =  \frac{\partial h_i}{\partial x} \cdot \frac{\partial F}{\partial u}$ for $i= 1,2$ with all partial derivatives evaluated $x = x_0$ and $u = 0$.  Recall that both $u$ and $|x-x_0|$ are assumed to be small; as such the higher order terms from \eqref{hopfcontrolpre} will be truncated to focus on the dominant contribution from the input.  A final  coordinate change $\hat{x} = K(\tilde{x} \cos(\theta_0) - \tilde{y} \sin(\theta_0))$, $\hat{y} = K(\tilde{x} \sin(\theta_0) + \tilde{y} \cos(\theta_0))$ will be considered.  Choosing $\theta_0$ so that $H_1 \sin(\theta_0) + H_2 \cos(\theta_0) = 0$, letting $K = 1/(H_1 \cos(\theta_0) - H_2 \sin(\theta_0))$, and truncating $ O(u^2)$ and $O(u |x - x_0|)$, the dynamics of the transformed coordinates become
\begin{align} \label{hopfcontrol}
    \dot{\hat{x}} &= \alpha \hat{x} - \beta \hat{y} + (a\hat{x} - b \hat{y} )(\hat{x}^2 + \hat{y} ^2) + u, \nonumber \\
    \dot{\hat{y} } &= \beta \hat{x} + \alpha \hat{y}  + (b\hat{x} + a \hat{y} )(\hat{x}^2 + \hat{y} ^2) ,
\end{align}
where ${\alpha}, {\beta}, {a},{b} \in \mathbb{R}$ are defined appropriately.  In the following analysis, Equation \eqref{hopfcontrol} will be referred to as the controlled Hopf normal form.

\subsection{Phase-Amplitude Reduction of the Controlled Hopf Normal Form}
While normal form theory guarantees the existence of a analytic coordinate changes that places the dynamics of \eqref{maineq} in the form \eqref{hopfcontrol} describing the system's behavior on the center manifold associated with the Hopf bifurcation, the specific coordinate the coefficients depend on the terms that comprise \eqref{maineq}.  This poses a challenge when $F$ is unknown.  In this case, the relevant terms must be inferred from data, but the mapping $y = g(x)$ is generally also unknown and is not guaranteed to be invertible.   Here, phase-amplitude reduction is considered to aid in model identification.  Considering the controlled Hopf normal form \eqref{hopfcontrol}, provided $a<0$ and $\alpha > 0$ the system has a stable periodic orbit $[\hat{x}^\gamma(\theta),\hat{y}^\gamma(\theta)] = [r_0 \cos(\theta),r_0 \sin(\theta)]$ where $r_0 = \sqrt{-\alpha/a}$.  Here, $\dot{\theta} = \omega$ when $u = 0$ where $\omega = \beta-\alpha b/a$.  Considering the adjoint equation from \eqref{adjeq} with the constraint that $\frac{d}{dt}[{\hat{x}},{\hat{y}}]\nabla \theta = \omega$ when evaluated on the periodic orbit taking $u = 0$, one can verify that 
\begin{equation} \label{phasehopf}
    \nabla \theta = \bigg[- \sqrt{-a/\alpha}\bigg( \sin(\theta) + b \cos(\theta)/a\bigg),  \sqrt{-a/\alpha} \bigg( \cos(\theta) - b \sin(\theta)/a \bigg)\bigg]^T,
\end{equation}
gives the gradient of the phase along the periodic orbit.  Towards considering the Floquet coordinates, first note that with the coordinate change $\hat{x} = r \cos(\theta), \hat{y} = r \cos(\theta)$, the dynamics of the radial coordinate are $\dot{r} = \alpha r + a r^3$.  For the fixed point at $r = r_0$, the associated linearization is $\dot{\Delta r} = -2 \alpha \Delta r$ giving the Floquet exponent 
\begin{equation} \label{kappaeq}
    \kappa_1 = -2 \alpha.
\end{equation}  
With this in mind, for the periodic orbit, $x^\gamma$, using Equation \eqref{isoeq2} one can verify that for the Floquet coordinate
\begin{equation} \label{amphopf}
    \nabla \psi_1 = C_1 [\cos(\theta), \sin(\theta)]^T,
\end{equation}
provides the gradient along the periodic orbit, where $C_1 \in \mathbb{C}$ is a constant that reflects the fact that the Floquet eigenfunctions are unique up to a constant scaling.  Considering the definition of $Z(\theta)$ and $I_1(\theta)$ provided after Equations \eqref{pred} and \eqref{isoeq}, respectively, in conjunction with the results from Equation \eqref{phasehopf} and \eqref{amphopf}, one finds
\begin{align} 
    Z(\theta) &= - \sqrt{-a/\alpha}\bigg( \sin(\theta) + b \cos(\theta)/a\bigg),  \label{zhopf} \\
    I_1(\theta) &= C_1 \cos(\theta),  \label{ihopf}
\end{align}
gives the phase and amplitude response curves for the transformed Hopf normal form from \eqref{hopfcontrol}.  It is important to mention that previous authors have considered the computation of the phase and amplitude response curve for the supercritical Hopf bifurcation \cite{brow04}, \cite{erme84}, \cite{cast13}, \cite{mong19}.  The specific representation used in Equations \eqref{zhopf} and \eqref{ihopf} will be used in the model identification strategy to follow.

\subsection{Preservation of the Terms of the Phase-Amplitude Reduction Following Transformation to the Controlled Hopf Normal Form} \label{shiftpreserve}

Consider the linearization near the periodic orbit \eqref{lineq}.   For any $x \in x^\gamma$ (so that $\psi_j = 0$ for all $j$), suppose $x$ is perturbed slightly to $x + \Delta x$ where $|\Delta x| = O(\epsilon)$ with $0<\epsilon \ll 1$.  To leading order in $\epsilon$, using the representation \eqref{floqform}, one can write
\begin{equation} \label{phasepert}
    x + \Delta x = x^\gamma(\theta) + \Delta \theta g_N + \sum_{j = 1}^{N-1} (\psi_j + \Delta \psi_j)g_j(\theta).
\end{equation}
Above, the terms $\Delta \theta$ and $\Delta \psi_j$ capture the change in the phase and Floquet coordinates in the basis of Floquet eigenfunctions and $g_N \equiv \frac{\partial x^\gamma}{\partial \theta}$ is the $N^{\rm th}$ Floquet eigenfunction of \eqref{maineq} which is associated with translation along the periodic orbit.  It is relatively straightforward to show that shifts in the phase and Floquet coordinate are preserved through the transformation to the controlled Hopf normal form.  To this end, let $D(x) = \begin{bmatrix} \hat{x}, \hat{y} \end{bmatrix}^T$ be the analytic, continuously differentiable coordinate transformation between the state variables in the given coordinate system and the controlled Hopf normal form.  In this transformed coordinate system, the periodic orbit becomes $\begin{bmatrix} \hat{x}^\gamma(\theta),\hat{y}^\gamma(\theta)\end{bmatrix}^T = D(x^\gamma(\theta))$.  To leading order accuracy, Equation \eqref{phasepert} becomes
\begin{align} \label{phasepert2}
\begin{bmatrix} \hat{x} + \Delta \hat{x} , \hat{y} + \Delta \hat{y} \end{bmatrix}^T &=  D \left(x^\gamma(\theta) +  \Delta \theta \frac{\partial x^\gamma}{\partial \theta} + \sum_{j = 1}^{N-1} \bigg[ (\psi_j + \Delta \psi_j)g_j(\theta) \bigg] \right)  \nonumber \\
&= D(x^\gamma(\theta)) + D^\prime(\theta) \cdot \frac{\partial x^\gamma}{\partial \theta} \Delta \theta +  \sum_{j = 1}^{N-1} \bigg[ (\psi_j + \Delta \psi_j) D^\prime(\theta) \cdot g_j(\theta) \bigg],
\end{align}
where $^\prime \equiv \frac{\partial}{\partial x}$.  Looking closer at the terms of \eqref{phasepert2}, one can rewrite Equation \eqref{phasepert2} as 
\begin{equation} \label{phasepert3}
    \begin{bmatrix} \hat{x} + \Delta \hat{x}, \hat{y} + \Delta \hat{y} \end{bmatrix}^T = D(x^\gamma(\theta)) + \hat{g}_N(\theta) \Delta \theta + \sum_{j = 1}^{N-1} \bigg[ C_j (\psi_j + \Delta \psi_j) \hat{g}_j(\theta) \bigg],
\end{equation}
where $\hat{g}_j$ for $j = 1, \dots, N$ are the Floquet eigenfunctions in the transformed coordinate system with the constant $C_j \in \mathbb{C}$ reflecting the fact that these eigenfunctions are unique to a constant scaling.  Above, the fact that $\hat{g}^N(\theta) \equiv \frac{\partial}{\partial \theta}\begin{bmatrix} \hat{x}^\gamma,\hat{y}^\gamma\end{bmatrix} = D^\prime(\theta) \cdot \frac{\partial x^\gamma}{\partial \theta}$ is used.  Comparing \eqref{phasepert} to \eqref{phasepert3}, one finds that phase shifts are preserved in the coordinate system associated with the controlled Hopf normal form and shifts in the Floquet coordinates are preserved up to a constant scaling.  Recalling that $\Delta x$ was arbitrary, this implies that $Z(\theta)$ is equivalent between coordinate systems with each $I_j(\theta)$ being preserved up to a constant scaling (i.e., the shape of $I_j(\theta)$ is the same).

\subsection{The Relationship Between the Geometry of Original and Transformed Coordinate Systems} \label{geomsec}

As shown here, the phase-amplitude reduction framework enables the inference of the terms associated with the controlled Hopf normal form despite the fact that the transformed coordinates $\hat{x}$ and $\hat{y}$ are generally not directly recoverable from measurements of the observable.  The key here is to relate the geometry of the problem relative to the coordinates of the controlled Hopf normal form from \eqref{hopfcontrol} to the geometry in observable space. Note that here and below, the asymptotic phase is denoted by $\hat{\theta}$ when viewed from the perspective of the controlled Hopf normal form and denoted by $\theta$ when viewed from the perspective of the original equations \eqref{maineq}.

\begin{figure}[htb]
\begin{center}
\includegraphics[height=2.2 in]{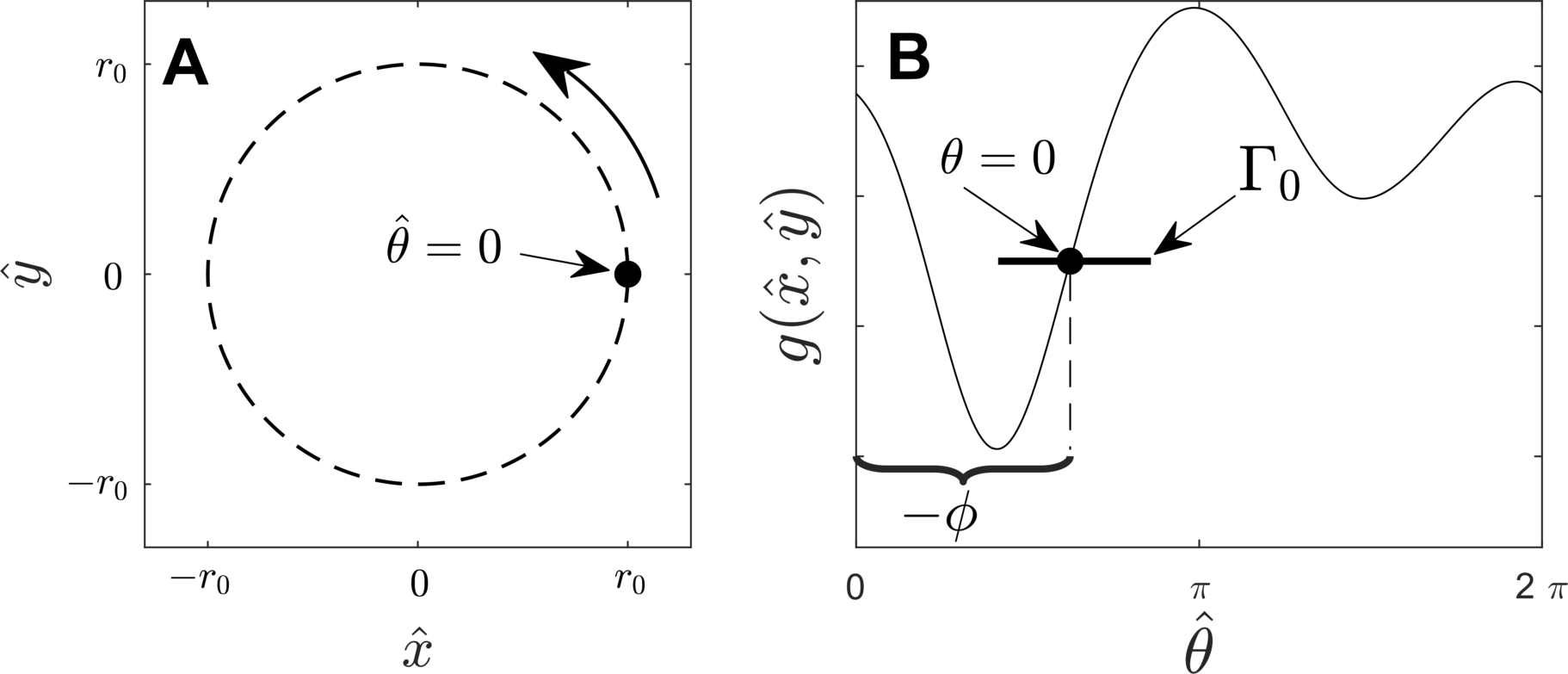}
\end{center}
\caption{Panel A highlights the geometry of the model identification problem from the perspective of the controlled Hopf normal form.  Panel B gives the representation as seen from the observable space in the context of the model identification strategies described in Section \ref{directback}.  Note that the the moment where $\theta = 0$ is determined by the choice of the Poincar\'e section $\Gamma_0$ and that $\theta = \hat{\theta} + \phi$.  This offset accounts for the fact that  $\hat{\theta} = 0$ will generally not correspond to  $\theta = 0$.}
\label{hopfgeometry}
\end{figure}

Figure \ref{hopfgeometry} highlights the details of the geometry in both frameworks for the case that $K = 1$ so that $y \in \mathbb{R}$.   Panel A emphasises the geometry of the Hopf normal form, with $(\hat{x},\hat{y}) = (0,r_0)$ corresponding to $\hat{\theta} = 0$.  From this perspective, the gradient of the phase and amplitude equations can be written analytically according to \eqref{phasehopf} and \eqref{amphopf}, respectively.  Panel B considers the output $g(\hat{x},\hat{y})$ plotted against $\hat{\theta}$ evaluated on the periodic orbit $[\hat{x}^\gamma(\theta),\hat{y}^\gamma(\theta)]$.  Considering the model identification strategies described in Section \ref{directback}, the terms of the phase-amplitude reduction can be obtained by first defining a Poincar\'e section $\Gamma_0$ transversal to $g(x^\gamma(\theta))$ at $\theta= 0$.  For practical purposes, $\Gamma_0$ is usually chosen to be an easily identifiable feature of the model output.  While the results from Section \ref{shiftpreserve} show that information about phase and Floquet coordinate shifts are preserved to a constant scaling,  there is no {\it a priori} way of identifying where $\hat{\theta} = 0$; as such a phase offset $\phi$ is necessary to relate $\theta$ and $\hat{\theta}$.  

\subsection{A Direct Method for Inference of the Unknown Coefficients of the Controlled Hopf Normal Form} \label{modelidsec}

The results presented above can be leveraged to develop a data-driven model identification strategy for systems that are close to a supercritical Hopf bifurcation.  This strategy can be implemented by applying as few as two pulse perturbations and examining the subsequent relaxation to the unforced periodic orbit.  As discussed in Section \ref{geomsec}, the notation $\theta$ and $\hat{\theta}$ is used to distinguish between the asymptotic phase in the original system coordinates and the coordinates of the controlled Hopf normal form, with the conversion $\theta = \hat{\theta} + \phi$.  The steps required to implement this model identification strategy are summarized below

\begin{enumerate}[1)]
\item Choose some $K-1$ dimensional hyperplane $\Gamma_0 \in \mathbb{R}^K$ transversal to $g(x^\gamma)$ to serve as a Poincar\'e section.  The output $g(x(t))$ crossing $\Gamma_0$ signifies that $\theta \approx 0$.

\item From passive observations (i.e.,~taking $u = 0$) let $T$ be defined as the average transversal between crossings of $\Gamma_0$.  Note that for noisy systems, the time between crossings of $\Gamma_0$ will generally not be constant.

\item After holding $u = 0$ long enough for transients to decay, at some time for which $g(x) \in \Gamma_0$ (so that $\theta \approx 0$), apply a short pulse $u(t) = m$ lasting $t_0$ time units.  By observing the subsequent relaxation to the periodic orbit, the strategies described in Section \ref{directback} can be used to provide an estimate of $\kappa_1$ (by observing the time between crossings of $\Gamma_0$ during relaxation to the periodic orbit) as well as $Z(0)$ and $\rho I_1(0)$ (according to Equations \eqref{phasedirect} and \eqref{floqdirect}, respectively).  Recall that $\rho$ is a constant that depends on system parameters that ultimately does not need to be inferred.

\item After the system fully relaxes to the periodic orbit, apply an identical pulse $T/4$ time units after the crossing of the $\Gamma_0$ Poincar\'e section (so that $\theta = \pi/2$.  Record the subsequent relaxation to the periodic orbit and use \eqref{phasedirect} and \eqref{floqdirect} to provide estimates of $Z(\pi/4)$ and $\rho I_1(\pi/4)$, respectively.  

\item  Considering Equation \eqref{ihopf}, 
\begin{equation} \label{ifit}
    I_1(\theta) = C_1 \cos(\theta - \phi).
\end{equation}
The two data points obtained from Steps 3 and 4 can be used to obtain an approximation for $\phi$ and $C_1/\rho$.  Note that due to symmetries in the cosine fitted function, the pair $(\phi,C_1/\rho)$ and $(\phi + \pi,-C_1/\rho)$ are equally good fits to the data; both of these options represent candidates for the true phase offset. 

\item Infer the terms $\alpha$, $\beta$, $a$, and $b$ for the controlled Hopf normal form starting with 
\begin{equation}
    \alpha = -\frac{\kappa_1}{2},
\end{equation}
which is obtained from \eqref{kappaeq}.   Considering Equation \eqref{zhopf}, one can write
\begin{equation} \label{zfit}
    Z(\theta) = - \sqrt{-a/\alpha}\bigg( \sin(\theta-\phi) + b \cos(\theta-\phi)/a\bigg).
\end{equation}
With knowledge of $\alpha$, and for the two candidates for $\phi$ from Step 5, a minimization over $a$ and $b$ can be performed with the data points obtained from from Steps 3 and 4.  Choose the true value of $a,b$, and $\phi$ to be the parameter set that minimizes the residual for the fit of $Z(\theta)$.  Consequently, the radius of the periodic orbit in controlled Hopf normal form coordinates is $r_0 = \sqrt{-\alpha/a}$.  Finally, considering the period on the unforced periodic orbit in relation to the radial form of the Hopf bifurcation, one finds 
\begin{equation}
    \beta = \frac{2 \pi}{T} - b r_0^2.
\end{equation}
\end{enumerate}
A few general notes about the implementation of the above model order reduction strategy and its subsequent use are provided below.

\begin{enumerate} [{Note 1)}]

\item Implementing Steps 1-6 above uniquely determines the coefficients associated with the controlled Hopf normal form.  For $x \in x^\gamma$, the output $g(x)$ crossing the $\Gamma_0$ level set corresponds to $(\hat{x}, \hat{y}) = (r_0 \cos(-\phi), r_0 \sin(-\phi))$.  Without further information, this is the only moment that an observation of the transformed coordinates can be obtained.  Under certain conditions, additional real-time estimates of the transformed state can be obtained as described in the section to follow.  

\item Only two pulse perturbations are sufficient to uniquely determine the unknown coefficients in the model fitting procedure described above. Pulses are applied at $\theta = 0$ and $\theta = \pi/2$ to ensure, in the worst case scenario, that $|I(\theta)|$ is within 70 percent of its peak value for at least one of the datapoints so that the dynamics associated with the amplitude coordinates can be  clearly observed.  Of course, additional data points can be incorporated into the parameter fitting from Steps 5 and 6.  Redundant measurements can be especially helpful to overcome measurement or system noise.  
\end{enumerate}

\subsection{Real-Time inference of State Variables for the Controlled Hopf Normal Form} \label{infsec}

In control applications, it can be useful to obtain real-time state information.  The state of the controlled Hopf normal form \eqref{hopfcontrol} cannot be directly observed from the system output $y$ as defined in Equation \eqref{maineq}.  For simplicity of exposition, it will be assumed that $K = 1$ so that $y \in \mathbb{R}$, but the following strategy could be straightforwardly applied in higher dimensions.  To proceed, noting that the series of near identity coordinate changes that yield $\tilde{x} = h_1(x)$ and $\tilde{y} = h_2(y)$ in the Hopf normal form \eqref{hnf} and the subsequent transformation to the controlled Hopf normal form  \eqref{hopfcontrol} are all invertible.  As such, letting $^{-1}$ denote the inverse, the function $D^{-1}(\hat{x},\hat{y}) = x$ which gives the state on the center manifold as a function of the coordinates in the controlled Hopf normal form is guaranteed to exist.  Assuming that both $\hat{x}$ and $\hat{y}$ are small and that $g$ is at least $C^1$ differentiable, to leading order one can write 
\begin{align} \label{yeq}
    y &= g(x) \nonumber \\
    &= g(D^{-1}(\hat{x},\hat{y})) \nonumber \\
    &\approx c_0 + c_1 \hat{x} + c_2 \hat{y},
\end{align}
where $c_0 = g(D^{-1}(0,0))$, $c_1 = \frac{\partial}{\partial \hat{x}} g(D^{-1})$, and $c_2 = \frac{\partial}{\partial \hat{y}} g(D^{-1})$, with partial derivatives taken at $\hat{x} = \hat{y} = 0$.  Taking the time derivative of \eqref{yeq} yields
\begin{align}
    \dot{y} &= c_1 \dot{\hat{x}} + c_2\dot{\hat{y}}  \nonumber \\
    &=  c_1 (\alpha \hat{x} - \beta \hat{y} + (a\hat{x} - b \hat{y} )(\hat{x}^2 + \hat{y} ^2) + u) + c_2 (\beta \hat{x} + \alpha \hat{y}  + (b\hat{x} + a \hat{y} )(\hat{x}^2 + \hat{y} ^2)) \nonumber \\
    &=  c_1 f_1(\hat{x},\hat{y}) + c_1 u + c_2 f_2(\hat{x},\hat{y}),
\end{align}
where the second line is obtained by substituting \eqref{hopfcontrol} and $f_1$ and $f_2$ in the third line are defined appropriately.  Still assuming $\hat{x}$ and $\hat{y}$ are small one can write
\begin{align} \label{ydoteq}
    \dot{y} &\approx \left( c_1 \frac{\partial f_1}{\partial \hat{x}} + c_2 \frac{\partial f_2}{\partial \hat{x}} \right) \hat{x} + \left(   c_1 \frac{\partial f_1}{\partial \hat{y}} + c_2 \frac{\partial f_2}{\partial \hat{y}} \right) \hat{y} + c_1 u \nonumber \\
    &= (c_1 \alpha + c_2 \beta) \hat{x} + (c_2 \alpha - c_1 \beta) \hat{y} + c_1 u \nonumber \\
    &= c_3 \hat{x} + c_4 \hat{y} + c_1 u,
\end{align}
where the second line is obtained by evaluating the partial derivatives at $\hat{x}=\hat{y} = 0$, noting that $f_1(0,0) = f_2(0,0)=0$ and defining $c_3$ and $c_4$ appropriately. Considering both \eqref{yeq} and \eqref{ydoteq} together, one can obtain an estimate of $\hat{x}$ and $\hat{y}$ according to
\begin{equation} \label{stateapprox}
    \begin{bmatrix}  \hat{x} \\ \hat{y} \end{bmatrix}  =  \begin{bmatrix} c_1 & c_2 \\ c_3 & c_4  \end{bmatrix}^{-1} \begin{bmatrix} y - c_0 \\ \dot{y} - c_1 u \end{bmatrix}.
\end{equation}
Above, the coefficients $c_0,c_1$ and $c_2$ can be obtained through observations of the periodic orbit $y(x^\gamma(\theta))$, comparing to $(\hat{x}(\theta), \hat{y}(\theta)) = (r_0 \cos(\theta-\phi), r_0 \sin(\theta-\phi))$, and inferring $c_0$, $c_1$, and $c_2$ through least squares fitting.  Instantaneous estimates of $\dot{y}$ can be  obtained from finite difference approximations.  

As mentioned previously, for $x \in x^\gamma$, when the output $g(x)$ crosses the $\Gamma_0$ level set, $(\hat{x}, \hat{y}) = (r_0 \cos(-\phi), r_0 \sin(-\phi))$ provided $\hat{x}$ and $\hat{y}$ are small.  The approximation \eqref{stateapprox} relies on the additional assumption that $g$ is smooth.  In practice, it can be useful to use \eqref{hopfcontrol} to approximate the state dynamics in response to the input $u$ and supplement with a state estimation algorithm with the information from \eqref{stateapprox}.  This will be described more carefully in the example from Section \ref{exampsec}.

\section{Examples} \label{exampsec}

\subsection{Model Identification for a Detailed Computational Model for Circadian Clock Oscillations}   \label{singexamp}
 As a preliminary example, a detailed model for mammalian circadian clock oscillations is considered.  The model contains 16 state variables that characterize the dynamical behavior of regulatory loops that govern the Per, Cry, Bmal1, and Clock genes \cite{lelo03}.  Full model equations are provided in Appendix \ref{apxa}.  The control parameter $v_{sP}(t) = v_{sP}^0 + u(t)$ is taken to be the control input which acts as a surrogate for the ambient light intensity. Here $v_{sP}^0 = 1.2 {\rm nM/h}$ corresponds to a moderate light intensity with higher and lower values corresponding to higher and lower intensity ambient light.  $M_P$ (the concentration of Per mRNA) is taken to be the measured output, i.e.,
 \begin{equation}
     y(t) = M_p(t).
 \end{equation}
Other observables can be used for this example and give comparable results.  This model undergoes a supercritical Hopf bufurcation when $u = -0.134$; when $u = 0$, the associated fixed point is unstable with eigenvalues $\lambda_{1,2} = 0.0254 \pm 0.275$.  This information is provided for context but is not used in the model identification strategy described below.

\begin{figure}[htb]
\begin{center}
\includegraphics[height=2.2 in]{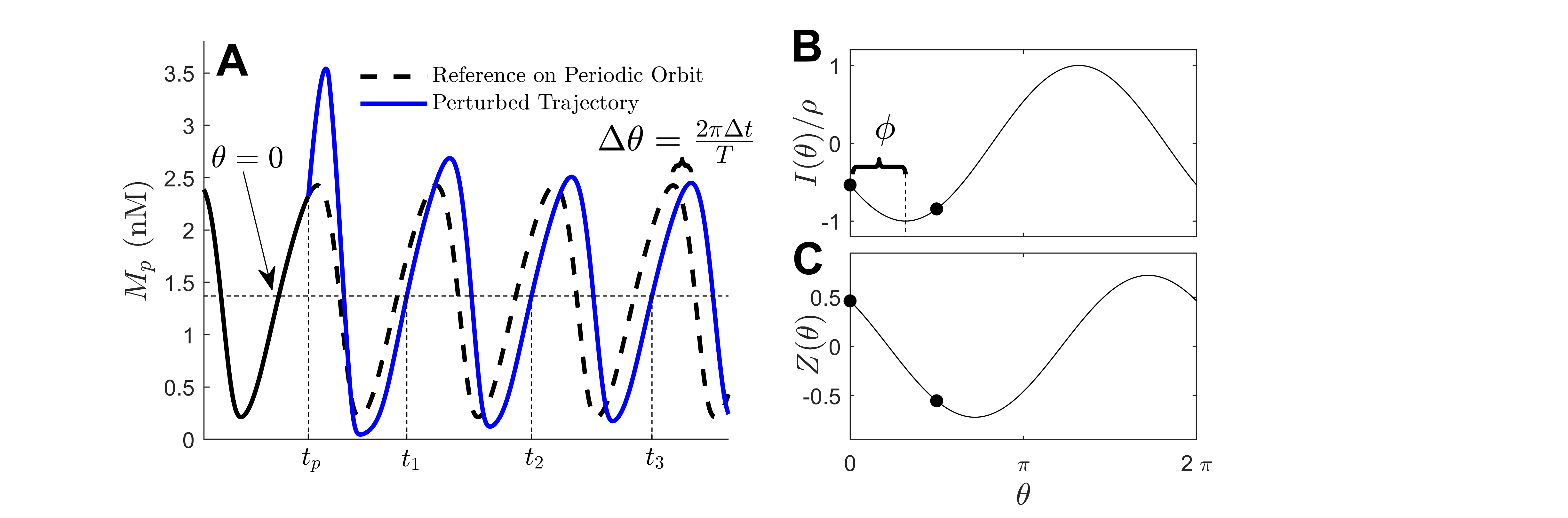}
\end{center}
\caption{The proposed model identification strategy is applied to the 16-dimensional circadian model described in Appendix \ref{apxa}.  In panel A, a pulse is applied at phase $\theta = \pi/2$ and the subsequent crossings of the $\Gamma_0$ surface are used to infer $\kappa_1$,  $I(\pi/2)/\rho$, and $Z(\pi/2)$. Direct estimates for $I(0)/\rho$ and $Z(0)$ are also obtained.  Black dots in panel B represent the direct estimates of $I(\theta)$ and the black line is a curve fit of the form \eqref{ifit}.  Panel C provides the same information for $Z(\theta)$ where the data is fit to a curve of the form \eqref{zfit}.  This information uniquely determines the unknown coefficients of the controlled Hopf normal form \eqref{hopfcontrol}.}
\label{showmodelid}
\end{figure}

Figure \ref{showmodelid} provides an illustration of the proposed model identification strategy applied to the circadian model from Appendix \ref{apxa}.  The $\Gamma_0$ Poincar\'e section is defined as the crossing of $M_p = 1.37$ with a positive slope.  $x^\gamma$ is taken to be the stable periodic orbit that results when taking $u = 0$, i.e.,~holding $v_{sP}$ constant at  $1.2 {\rm nM/h}$.  In panel A, a short pulse in $u$ is applied at $\theta = \pi/2$, $T/4$ hours after the phase crosses the $\Gamma_0$ surface.  Subsequent crossings of the $\Gamma_0$ surface are used to infer $\kappa_1$, $I(\pi/2)/\rho$ and $Z(\pi/2)$ using the strategy described in Section \ref{directback}.  A second pulse applied at $\theta = 0$ is used to infer $I(0)/\rho$ and $Z(0)$.  Recall that while the data driven model identification strategy from Section \ref{directback} provides an estimate of $\rho I(\theta)$ where $\rho$ is an undetermined constant.  The two directly inferred data points shown as black lines in panel B (resp.,~C) are used to fit curves of the form \eqref{ifit} (resp.,~\eqref{zfit}).   These fits are used to determine the unknown coefficients $a = -0.0106$, $b = -0.0034$, $\alpha = 0.0224$, and $\beta = 0.2721$ in the controlled Hopf normal form \eqref{hopfcontrol} as well as the phase offset $\phi = 1.005$ rad.  Note here that the model dynamics are inferred by observing the output during the recovery to the periodic orbit for only two pulse perturbations.

The resulting model is validated with the formulation and solution of an optimal control problem of shifting the oscillation timing, with relevance to the development light exposure scheduling for jet-lag mitigation strategies \cite{maye23}, \cite{chri20}, \cite{serk14}.  Following the optimal control problem formulation described in Appendix \ref{apxb} and letting $\xi(t_i) = \begin{bmatrix} \hat{x}(t_i),\hat{y}(t_i) \end{bmatrix}^T$, the cost function considered here is 
\begin{equation} \label{costcirc16}
    J = k \sum_{i = 0}^\eta \left[ 1- \exp( -30  || \xi(i \Delta t) - \xi_{\rm targ}(i \Delta t)||^2)   \right]   +   \sum_{i = 0}^{\eta-1}  u^2 (i \Delta t),
\end{equation}
where $k$ is a positive constant that weights the relative importance of the  state-based and control-based costs and $\xi_{\rm targ}(t) = r_0 \begin{bmatrix} \cos(\theta_0 + \omega t, \sin(\theta_0 + \omega t)  \end{bmatrix}$ with $\omega = \beta-\alpha b/a$ and $r_0 = \sqrt{-\alpha/a}$ and $\theta_0 \in [0,2\pi)$.  The target $\xi_{\rm targ}$ represents a state evolving on the periodic orbit at the unperturbed natural frequency.  During the application of the optimal control, the state estimate is updated according to
\begin{equation} \label{stateest}
    \zeta_{i+1} = f(\zeta_i,u) + \nu \left(    \begin{bmatrix} c_1 & c_2 \\ c_3 & c_4  \end{bmatrix}^{-1} \begin{bmatrix} y - c_0  \\ \dot{y} - c_1 u \end{bmatrix} - f(\zeta_i,u) \right).
\end{equation}
Above, as in $\eqref{zetaeq}$, $\zeta_i = [ \hat{x}((i \Delta t),\hat{y}(i \Delta t)]^T$ and $f$ gives the evolution of $\hat{x}$ and $\hat{y}$ under the evolution of \eqref{hopfcontrol} when applying a constant $u$.  As described in Section \ref{infsec}, the constants $c_1 - c_4$ can be obtained through least-squares fitting with knowledge of $y(x^\gamma(\theta))$.  The constant $\nu \in (0,1]$ sets how aggressively to update the state when comparing the expected evolution of the coordinates in the controlled Hopf normal to estimates obtained from  direct observations of $y$ and $\dot{y}$.  For this example, $\nu = 0.02$.

\begin{figure}[htb]
\begin{center}
\includegraphics[height=2.2 in]{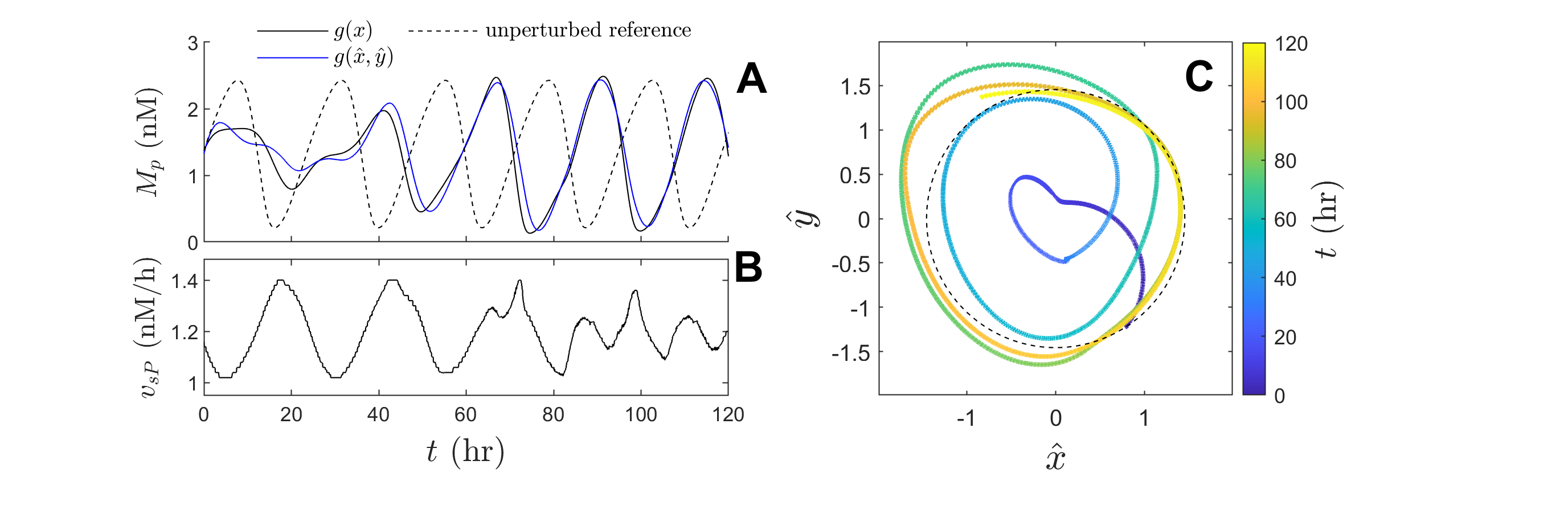}
\end{center}
\caption{After model identification, the optimal control strategy from Appendix \ref{apxb} is applied to the inferred model \eqref{hopfcontrol}.  The goal is to shift the phase 12 hours in advance while limiting the control effort.  The optimal control is shown in panel B.  The solid black line in panel A gives the resulting model output with the blue line obtained using the estimation of the state of the inferred model \eqref{hopfcontrol}. The dashed line gives the trajectory for a system for which $u = 0$ for reference.  Panel C shows the evolution of $\hat{x}$ and $\hat{y}$ over time, estimated according to \eqref{stateest},  during the application of the control.  The dashed circle shows the unperturbed periodic orbit with radius $r_0 = \sqrt{-\alpha /a}$ for reference.}
\label{circphaseshift}
\end{figure}

As described in Appendix \ref{apxb}, the cost-to-go function \eqref{ctg} associated with the cost function \eqref{costcirc16} is solved numerically taking $\Delta t = 0.1$ hours, $\eta  =1200$, $k = 1/50$, and $\theta_0 = -\phi + 12 \omega$ with allowable control input $u \in [-0.2, 0.2]$.  For an initial condition at $\theta = 0$ on $x^\gamma$ (i.e.,~the periodic orbit that results when taking $u = 0$)  the resulting control defined by \eqref{optimalinput} will advance the state by 12 hours, or approximately $\pi$ radians.  Figure \ref{circphaseshift} shows the result of the application of this optimal control strategy.  The black trace in Panel A gives the actual output for the model with the blue trace giving the output computed according to \eqref{yeq} as a function of the coordinates of the controlled Hopf normal form.  Note that these two traces are slightly different and that the correction term from the state update equation \eqref{stateest} is small.  The dashed line gives a plot of the unperturbed reference highlighting the shift in the oscillation timing caused by the optimal input which is shown in Panel B.  After the application of the optimal input, the true shift in oscillation timing as determined by successive crossings of the $\Gamma_0$ Poincar\'e surface is 11.3 hours.  The curve in Panel C shows the evolution of $\hat{x}$ and $\hat{y}$ over time.  Note that these coordinates cannot be measured directly and are instead estimated according to \eqref{stateest}.  For the large phase shift, the state comes close to the unstable fixed point at $\hat{x} = \hat{y} = 0$ before being driven back to the periodic orbit with the correct phase.

\FloatBarrier

\subsection{Model Identification for Population-Level Oscillations in a Coupled Oscillator Model}  \label{circcontrolsec}
As a second example, a phenomenological model is considered for coupled oscillations in the suprachiasmatic nucleus (SCN) \cite{gonz05}, the ``master clock" responsible for circadian pacemaking:
\begin{align} \label{popmodel}
\dot{B_i} &= v_1 \frac{K_1^n}{K_1^n + D_i^n} - v_2 \frac{B_i}{K_2+B_i}  +     h_c \frac{K F}{K_c + K F}  +  \sigma_i u(t) + \sqrt{2 D} \eta_i(t),\nonumber \\
\dot{C_i} &= k_3 B_i - v_4 \frac{C_i}{K_4+C_i}, \nonumber \\
\dot{D_i} &= k_5 C_i - v_6 \frac{D_i}{K_6+D_i}, \nonumber \\
\dot{E_i} & = k_7 B_i - v_8 \frac{E_i}{K_8 + E_i}, \quad i = 1,\dots,N.
\end{align}
In mammals, the SCN is comprised of roughly 10,000 coupled neurons that respond to environmental cues such as light intensity in order to entrain to a 24-hour light-dark cycle \cite{repp02}, \cite{moor02}.  The above model considers $N = 3000$ oscillators with the variables $B_i$, $C_i$, and $D_i$ representing mRNA concentrations of clock gene, the associated protein, and the nuclear form of the protein, respectively, for the $i^{\rm th}$ oscillator. $E_i$ represents a neurotransmitter that contributes to the mean-field coupling $F(t)  = \frac{1}{N} \sum_{i = 1}^N E_i(t)$ that influences the behavior of the variable $B$.  $F(t)$ is taken to be the measured output for this model.  Time is measured in units of hours.   The control input $u$ can be thought of as the intensity of ambient light; each neuron has an intrinsic sensitivity to light $\sigma_i = \max(1 + 0.4N(0,1),0)$ where $N(0,1)$ is a normal distribution with zero mean and unit variance.  The term $\sqrt{2 D} \eta_i(t)$ incorporates an independent and identically distributed zero-mean white noise process with intensity $D = 0.0001$.  Nominal model parameters are $n = 5, v_1 = 0.55, v_2 = 0.39, v_4 = .35, v_6 = .35, v_8 = 1, k_3 = .7, k_5 = .5,  k_7 = .35, K_1 = 1, K-2 = 1, K_4 = 1, K_6 = 1, K_8 = 1, h_c = .35, K_c = 1$ and $K = .5$. To incorporate heterogeneity, parameters $v_1, v_2, v_4, v_6, k_3$, and $k_5$ are drawn from a normal distribution with mean being equal to the nominal parameter and a variance of 0.0001.  

\begin{figure}[htb]
\begin{center}
\includegraphics[height=1.8 in]{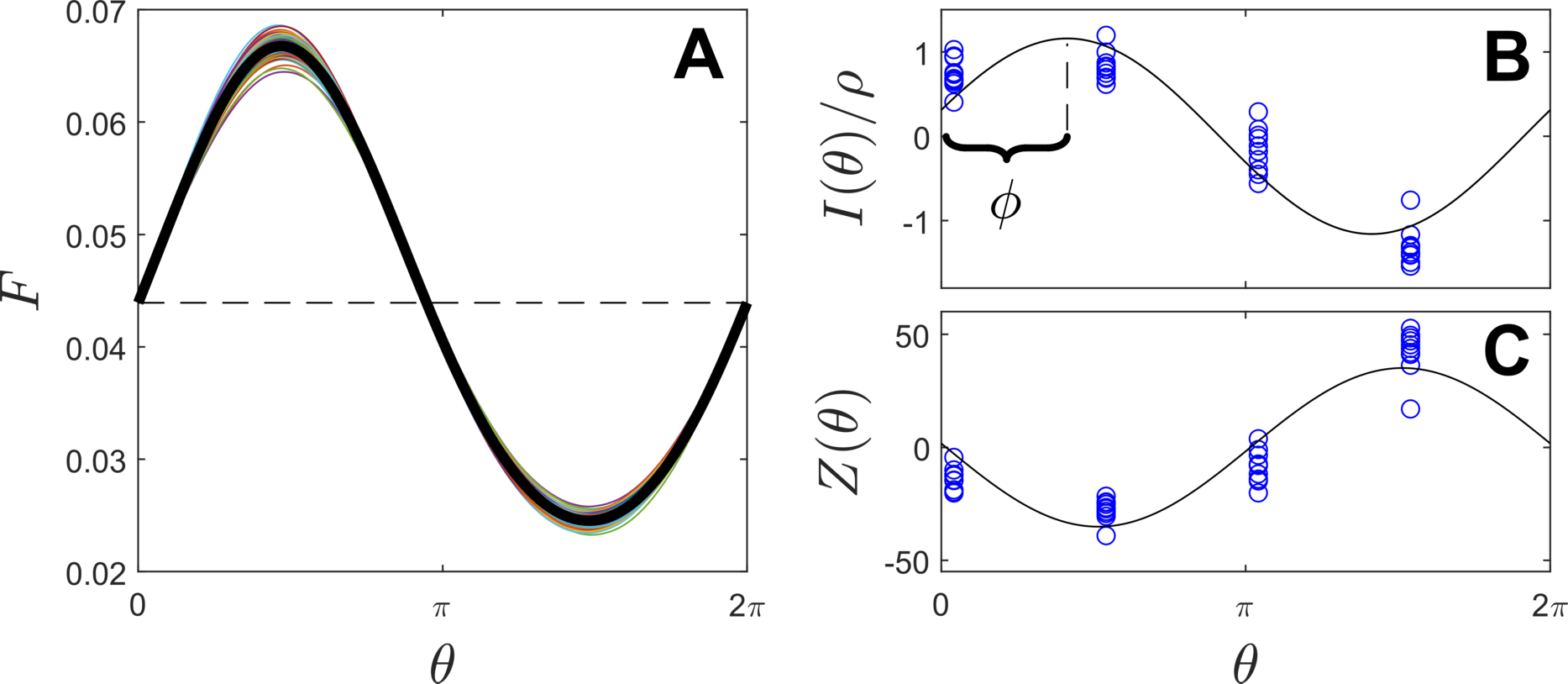}
\end{center}
\caption{For the population circadian model \eqref{popmodel}, colored traces in panel A show recordings of the output, $F(t)$, between crossings of the $\Gamma_0$ Poincar\'e section.  $F(x^\gamma)$ is taken to be the average of these traces.  Pulse inputs are applied with resulting estimates of $I(\theta)/\rho$ and $Z(\theta)$ shown as blue circles in panels B and C, respectively.  This data is used to fit sinusoidal curves shown in black which in turn are used to infer the terms of the controlled Hopf normal form from \eqref{hopfcontrol}.}
\label{circpopdata}
\end{figure}

Figure \ref{circpopdata} provides an illustration of the proposed model identification strategy applied to the circadian model from \eqref{popmodel}.  The $\Gamma_0$ Poincar\'e section is taken to correspond to the crossing of $F(t) = 0.044$ with a positive slope and is represented by the dashed line in panel A.  Recording 230 cycles, the mean return time to the $\Gamma_0$ section is 24.17 hours.  Colored traces in panel A show individual recordings of $F$ over these individual cycles. $F(x^\gamma)$ is taken to be the average of these 230 recordings and shown as a solid black line.  Note that because of the noise, the system never settles to the periodic orbit.   Following the strategy from Section \ref{directback}, pointwise estimates of $I(\theta)$ and $Z(\theta)$ are obtained applying a short pulse of input at specific times and recording the relaxation back to the limit cycle.  Results of individual trials in panels B and C are shown with blue circles and the black curves are fit to functions of the form \eqref{ifit} and \eqref{zfit} for $I(\theta)/\rho$ and $Z(\theta)$, respectively (recall from Section \ref{directback} that $I(\theta)$ can only be estimated up to a positive constant $\rho$).  Here, multiple samples are considered in order to contend with the noise inherent to \eqref{popmodel} but similar fits to the data can be obtained using only a subset of these samples. The fitted curves from panels B and C are used to infer the terms of the controlled Hopf normal form from \eqref{hopfcontrol}; here $a = -0.41$, $b = -1.24$, $\alpha = 0.0034$, and $\beta = 0.27$ with phase offset $\phi = 1.30$.  

\begin{figure}[htb]
\begin{center}
\includegraphics[height=1.9 in]{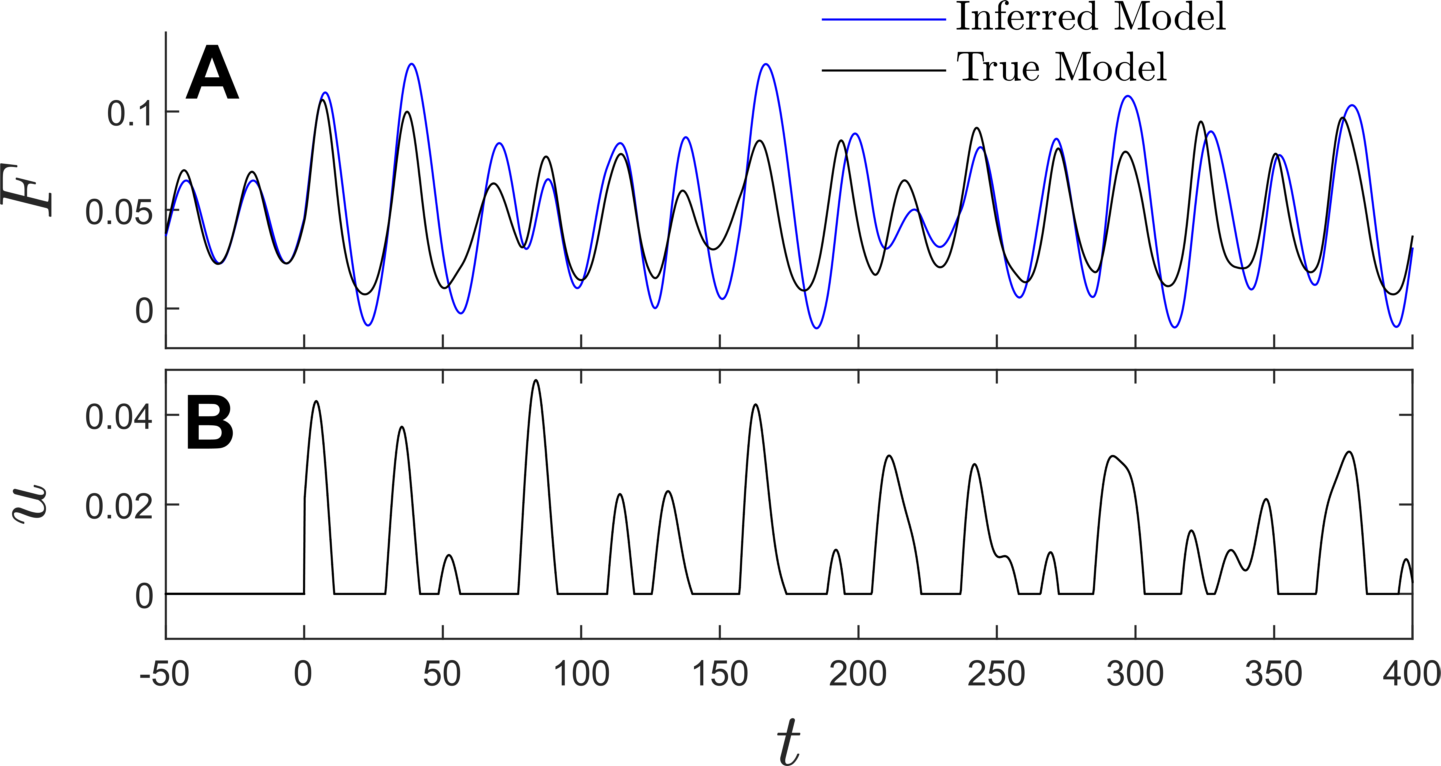}
\end{center}
\caption{Panel A shows the output from both the true model and inferred model in response the applied input from panel B.  While not a perfect match, the 2-dimensional inferred model provides a reasonable representation of the dynamics of the 12,000 dimensional model \eqref{popmodel}.}
\label{poptest}
\end{figure}

The accuracy of the inferred model resulting from the proposed model identification procedure is first tested by applying an the open loop input $u(t) = \max [ .02 \sin(0.24 t) + .02 \cos(0.15 t) + .025 \sin(0.4 t),0]$ starting at $t = 0$.  Results are shown in Figure \ref{poptest}.  While the output and input do not match perfectly, the timing of the oscillations in response to the input matches well.  Results are qualitatively similar when using other inputs of similar magnitude.  

\begin{figure}[htb]
\begin{center}
\includegraphics[height=2.3 in]{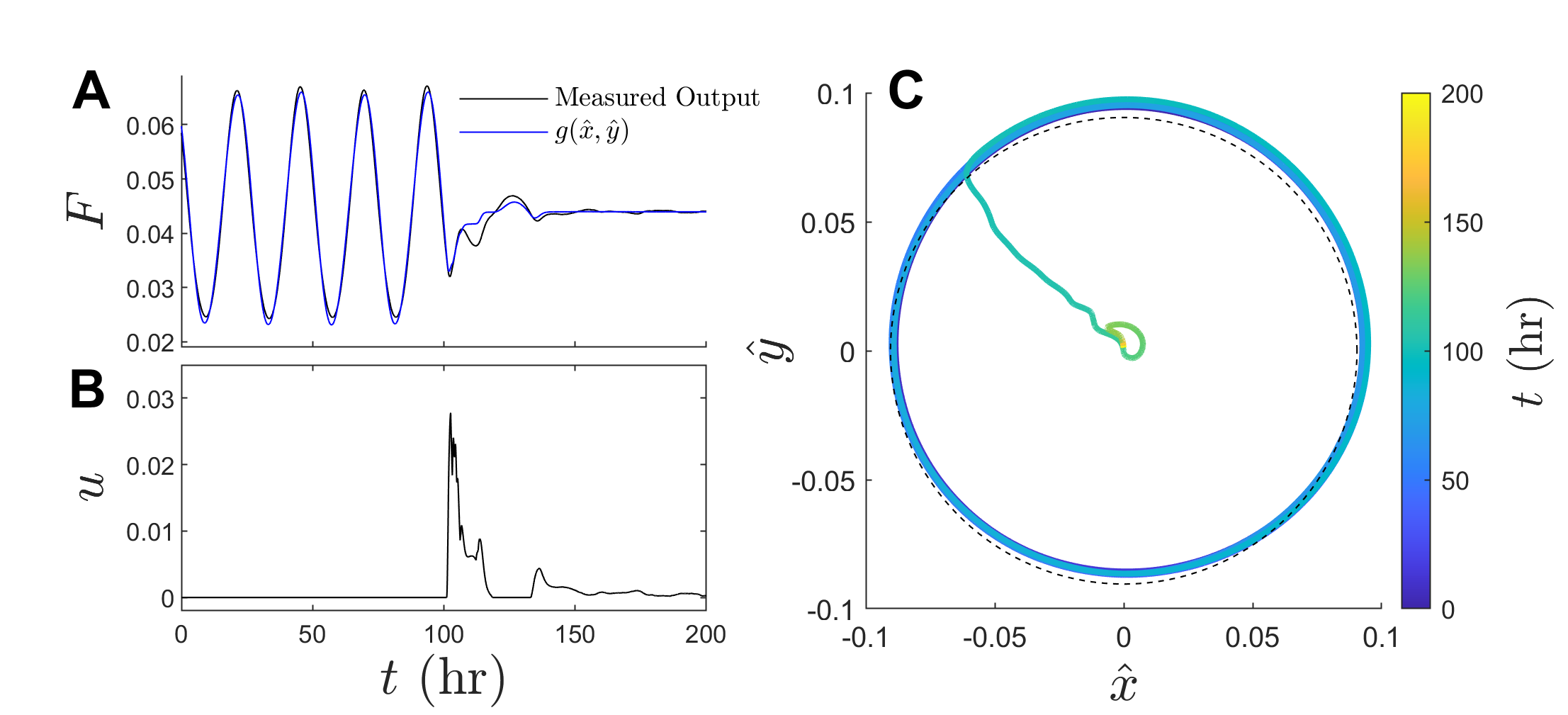}
\end{center}
\caption{Optimal control applied to the population model for circadian oscillations \eqref{popmodel}.  The optimal state feedback controller is computed according to the dynamic programming strategy described in the text. The goal here is to drive the system to the phaseless set, i.e.,~the unstable fixed point associated with the Hopf bifurcation.  At $t = 100$ the controller is turned on.  Panel A shows the output from the full model (black line) along with an estimate obtained from \eqref{yeq} (blue line) which uses the estimate for the state of the controlled Hopf normal form.  Panel B shows the control applied as a function of time.  Panel C shows the evolution of $\hat{x}$ versus $\hat{y}$ during the simulation.  The  unperturbed periodic orbit with radius $r_0 = \sqrt{-a/\alpha}$ is shown for reference as a dashed line.}  
\label{circpopcont}
\end{figure}

The inferred model is also validated by formulating and solving an optimal control problem with the goal of driving the system to the unstable fixed point of the associated Hopf bifurcation (i.e.,~the phaseless set \cite{winf01}).  Following the optimal control formulation described in Appendix \ref{apxb} and letting and letting $\xi(t_i) = \begin{bmatrix} \hat{x}(t_i),\hat{y}(t_i) \end{bmatrix}^T$, the cost function considered here is 
\begin{equation} \label{circpopcost}
    J = k \sum_{i = 0}^\eta \left[ 1-\exp(-20 || \xi(i \Delta t) ||^2) \right] +  \sum_{i = 0}^{\eta-1} u^2(i \Delta t),
\end{equation}
where $k$ is a positive constant that weights the relative importance of the state-based and control-based costs. As described in Appendix \ref{apxb} the cost-to-go function \eqref{ctg} associated with the cost function \eqref{circpopcost} is solved numerically taking $\Delta t = 0.1$ hours, $\eta = 240$, and $k = 1$ with allowable control input $u \in [0,0.04]$.   The resulting optimal control is applied to the full order model \eqref{popmodel} with a 24 hour prediction horizon, i.e.,~by using the $J_1$ cost-to-go function at all times when computing the optimal control according to \eqref{optimalinput}.  During the application of the optimal control, the inferred state of the controlled Hopf normal form is updated according to \eqref{stateest} where $\nu = 0.005$.   Results are shown in Figure \ref{circpopcont}.  Starting from an initial condition near the unperturbed periodic orbit, the control is turned on at $t = 100$ hours and quickly drives the system to the phaseless set.  Panel A shows the controlled output (black line) for the model along with the output computed according to \eqref{yeq} as a function of the coordinates of the controlled Hopf normal form (blue line); notice that there is relatively little discrepancy between these two curves.  The applied optimal control is shown in panel B.  Panel C shows the evolution of $\hat{x}$ and $\hat{y}$ over time calculated according to the estimator \eqref{stateest}.  The initial burst of input at $t = 100$ quickly drives the system close to the unstable fixed point of the controlled Hopf normal form.  The subsequent inputs are much smaller, counteracting the effect of noise to keep the system close to the target.  

\FloatBarrier

\section{Discussion and Conclusion} \label{concsec}

In this work, a general strategy is developed and investigated for data-driven inference of oscillatory dynamical systems close to a supercritical Hopf bifurcation.  By first augmenting the Hopf normal form with a control input to yield \eqref{hopfcontrol}, the unknown coefficients can be determined uniquely by obtaining as few as two pointwise estimates of the phase and amplitude response curves when considering the dynamics from the perspective of isochrons and Floquet coordinates.  Examples illustrate the utility of this model identification strategy with the formulation and solution of different optimal control problems.  

With regard to the underlying form of the model from \eqref{maineq}, while the proposed model identification strategy assumes that the underlying state dynamics that comprise $F(x,u)$ are sufficiently smooth, there is no explicit requirement about differentiability or continuity of the state-to-output map $g(x)$.  In principle, sparse identification strategies \cite{brun16b}, \cite{pant19}, \cite{fase21} could also be used for data-driven model identification purpose, but such strategies may prove difficult in situations where $g$ is not sufficiently smooth.  Note that the real-time state inference strategy from Section \ref{infsec} does require that $g$ is at least $C^1$ differentiable, but even without the state-inference strategy, the phase and associated coordinates of the controlled Hopf normal form can always be estimated once per cycle at the crossing of the $\Gamma_0$ Poincar\'e section provided the system is close enough to the unperturbed periodic orbit.

There are many limitations of the proposed model identification strategy that are worth mentioning.  The proposed model identification approach cannot be straightforwardly generalized to systems undergoing a supercritical Hopf bifurcation.  For such systems, the resulting stable limit cycle typically settles far from the fixed point associated with the Hopf bifurcation thereby invalidating the assumptions required to implement the proposed approach.  For these systems, other data-driven approaches that rely on related phase-based model identification approaches would likely be more applicable \cite{wils23direct}, \cite{wils20ddred}, \cite{wils18operat}.   Additionally, the current manuscript does not explicitly consider noise or other sources of uncertainty in the problem formulation.  In the results from Section \ref{circcontrolsec} noise is included in the model equations necessitating multiple redundant measurements of the associated phase and amplitude response curves with the goal of averaging out the effects of noise.  It would be of particular interest to obtain rigorous estimates of how noise contributes to errors in the representation of the resulting models.  Furthermore, while the proposed approach is valid provided the system is sufficiently close to a subcritical Hopf bifurcation, estimates for how this estimation degrades as the state moves farther from the Hopf bifurcation are not considered.  

Almost certainly, other data-driven model identification strategies could be applied to the systems considered in this work, but careful consideration of the dynamical behavior near the Hopf bifurcation yields a relatively simple model fitting approach that requires very little data.  The results presented here serve to highlight the need for the continued development of data-driven model identification strategies that leverage the knowledge of normal forms near bifurcations of vector fields.  While the strategy proposed in this work is only applicable for one specific (but commonly observed) bifurcation, other authors have investigated the shape of phase and amplitude response curves for oscillatory systems that are near other bifurcations \cite{erme96},  \cite{brow04}, \cite{mong21}.  It would be of interest to consider the development of similar model identification approaches for systems near other types of bifurcations.

This material is based upon the work supported by the National Science Foundation (NSF) under Grant No.~CMMI-2140527.


\begin{appendices}

\section{Circadian Model Equations} \label{apxa}
\renewcommand{\thetable}{A\arabic{table}}  
\renewcommand{\thefigure}{A\arabic{figure}} 
\renewcommand{\theequation}{A\arabic{equation}} 
\setcounter{equation}{0}
\setcounter{figure}{0}

The circadian oscillator model used in Section \ref{singexamp} was originally published in \cite{lelo03}.  The version used here is comprised of 16 coupled ordinary differential equation.  The state variables as follows:~Concentrations of  Per, Cry, and Bmal1 mRNA are designated by $M_P$, $M_C$, and $M_B$, respectively;  phosphorylated (resp.,~nonphosphorylated)  Per and Cry proteins in cytosol are designated by $P_{CP}$ and $C_{CP}$ (resp.,~$P_C$ and $C_C$);  concentrations of Per-Cry complex in cytosol and nucleus are designated by $PC_C$, $PC_N$, $PC_{CP}$, and $PC_{NP}$;  concentrations of BMAL1 in cytosol and nucleus are designated by $B_C$, $B_{CP}$, $B_N$, and $B_{NP}$;  the variable $I_N$ represents the inactive complex between Per-Cry and Clock-Bmal1 in the nucleus.  Subscripts $C$, $N$, $CP$ and $NP$ denote cytosolic, nuclear, cytosolic phosphorylated, and phosphorylated forms, respectively.    The model equations are:
\begin{align}
\dot{M}_P &= v_{sP}(t) \frac{B_N^n}{K_{AP}^n + B_N^n} - v_{mP} \frac{M_P}{K_{mP} + M_P} - k_{dmp}M_P, \label{a1} \\
\dot{M}_C &= v_{sC}\frac{B_N^n}{K^n_{AC} + B_N^n} - v_{mC} \frac{M_C}{K_{mC} + M_C} - k_{dmc} M_C, \\
\dot{M}_B &= v_{sB}\frac{K_{IB}^m}{K_{IB}^m + B_N^m} - v_{mB} \frac{M_B}{K_{mB} + M_B} - k_{dmb} M_B, \\
\dot{P}_C &= k_{sP}M_P - V_{1P} \frac{P_C}{K_p + P_C} + V_{2P}\frac{P_{CP}}{K_{dp} + P_{CP}} + k_4 PC_C - k_3 P_CC_C - k_{dn}P_C,\\
\dot{C}_C &= k_{sC}M_C - V_{1C} \frac{C_C}{K_p + C_C} + V_{2C}\frac{C_{CP}}{K_{dp} + C_{CP}} + k_4 PC_C - k_3 P_CC_C - k_{dnc}C_C,\\
\dot{P}_{CP} &= V_{1P} \frac{P_C}{K_p + P_C} - V_{2P} \frac{P_{CP}}{K_{dp} + P_{CP}} - v_{dPC} \frac{P_{CP}}{K_d + P_{CP}} - k_{dn}P_{CP},\\
\dot{C}_{CP} &= V_{1C} \frac{C_C}{K_p + C_C} - V_{2C} \frac{C_{CP}}{K_{dp} + C_{CP}} - v_{dCC} \frac{C_{CP}}{K_d + C_{CP}} - k_{dn}C_{CP},\\
\dot{PC}_C &= - V_{1PC} \frac{PC_C}{K_p + PC_C} + V_{2PC}\frac{PC_{CP}}{K_{dp} + PC_{CP}} - k_4 PC_C + k_3 P_C C_C \nonumber \\
&\quad+ k_2 PC_N - k_1PC_C - k_{dn}PC_C,\\
\dot{PC}_N &= - V_{3PC} \frac{PC_N}{K_p + PC_N} + V_{4PC}\frac{PC_{NP}}{K_{dp} + PC_{NP}} - k_2 PC_N + k_1 PC_C \nonumber \\
&\quad- k_7B_N PC_N + k_8I_N - k_{dn}PC_N,\\
\dot{PC}_{CP} &= V_{1PC} \frac{PC_C}{K_p + PC_C} - V_{2PC} \frac{PC_{CP}}{K_{dp} + PC_{CP}} - v_{dPCC} \frac{PC_{CP}}{K_d + PC_{CP}} - k_{dn}PC_{CP},\\
\dot{PC}_{NP} &=- V_{3PC} \frac{PC_N}{K_p + PC_N} - V_{4PC}\frac{PC_{NP}}{K_{dp} + PC_{NP}} - v_{dPCN} \frac{PC_{NP}}{K_d + PC_{NP}} - k_{dn}PC_{NP},  \\
\dot{B}_C &= k_{sB} M_B - V_{1B} \frac{B_C}{K_p + B_C} + V_{2B} \frac{B_{CP}}{K_{dp} + B_{CP}} - k_5B_C + k_6B_N - k_{dn}B_C,\\
\dot{B}_{CP} &=  V_{1B} \frac{B_C}{K_p + B_C} - V_{2B} \frac{B_{CP}}{K_{dp} + B_{CP}} - v_{dBC} \frac{B_{CP}}{K_d + B_{CP}} - k_{dn}B_{CP},\\
\dot{B}_{N} &=  -V_{3B} \frac{B_N}{K_p + B_N} + V_{4B} \frac{B_{NP}}{K_{dp} + B_{NP}} + k_5 B_C - k_6 B_N - k_7 B_N PC_N  \nonumber \\
& \quad + k_8 I_N - k_{dn}B_N, \\
\dot{B}_{NP} &=  V_{3B} \frac{B_N}{K_p + B_N} - V_{4B} \frac{B_{NP}}{K_{dp} + B_{NP}} - v_{dBN} \frac{B_{NP}}{K_d + B_{NP}} - k_{dn}B_{NP},\\
\dot{I}_N &= -k_8 I_N + k_7 B_N PC_N-v_{dIN} \frac{I_N}{K_d + I_N} - k_{dn}I_N  \label{a16}.
\end{align}
Basal values listed in Supplementary Table 1 of \cite{lelo03} are used with the exception of $k_1 = 0.58$ and $k_2 = 2.0$, which determine the dynamics of the nonphysphorylated cytosolic Per and Cry protein concentrations.  Units of time are in hours.  Light intensity acts as an input for this model that changes the value of $v_{sP}$, the maximum rate of Per expression.  $v_{sP} = 1.2$ nM/h corresponds to a moderate light intensity with higher/lower values corresponding to higher/lower light intensities.

\section{Optimal Control Using a Dynamic Programming Approach} \label{apxb}
\renewcommand{\thetable}{B\arabic{table}}  
\renewcommand{\thefigure}{B\arabic{figure}} 
\renewcommand{\theequation}{B\arabic{equation}} 
\setcounter{equation}{0}
\setcounter{figure}{0}

The inferred models from Sections \ref{singexamp} and \ref{circcontrolsec} are validated by formulating and solving an optimal control problem. In this work, a dynamic programming approach is used.  Note that the models obtained from the model identification procedure in Section \ref{modelidsec} capture the behavior in response to general inputs and that if desired, other control strategies could be considered.  

To begin, the inferred model from \eqref{hopfcontrol} is discretized considering time steps of length $\Delta t$.  For some initial condition $\zeta_0 = [\hat{x}(0), \hat{y}(0)]^T $, define $\zeta_i = [\hat{x}(i \Delta t), \hat{y}(i \Delta t)]^T$ to be the solution to \eqref{hopfcontrol}.  Taking $u$ to be constant on the interval $t = [i \Delta t, (i+1) \Delta t]$ one can use \eqref{hopfcontrol} to numerically evaluate the mapping 
\begin{equation}\label{zetaeq}
    \zeta_{i+1} = f(\zeta_i,u).
\end{equation}
Next, one can define a general cost function
\begin{equation}
J \bigg( (u_i)_{i = 0}^{\eta-1},  (\zeta_i)_{i = 0}^\eta \bigg) = \sum_{i = 0}^{\eta-1} c_1 (u_i,i) + \sum_{i = 0}^\eta c_2(\zeta_i,i),
\end{equation}
where $c_1$ and $c_2$ are time-dependent functions of the input and state, respectively.  The optimization problem considered in this work is defined as follows:~for a sequence of inputs $u_0,u_1,u_{\eta-1}$ where the input $u_i$ denotes a constant input applied over the interval $[i \Delta t,(i+1)\Delta t)]$
\begin{equation}
\underset{\substack{u_i \\ 0 \leq i \leq \eta-1 }}{\rm minimize} \;  J \bigg( (u_i)_{i = 0}^{\eta-1},  (\zeta_i)_{i = 0}^\eta \bigg),
\end{equation}
over all allowable inputs $u_i \subset \mathcal{U}$ where $\eta$ is the number of iterations considered and the evolution of $\zeta$ is giverned by \eqref{zetaeq}.  As described in \cite{bell66}, \cite{kirk98}, The minimization can be performed using a dynamic programming approach.  Defining
\begin{equation} \label{ctg}
J^*_{\eta-y} (\zeta_{\eta-y} )= \begin{cases}  \min\limits_{  \substack {u_i  \\  \eta-y \leq i \leq \eta-1  }}    J \bigg( (u_i)_{i = \eta-y}^{\eta-1},  (\zeta_i)_{i = \eta-y}^\eta \bigg), & \text{if } y>0,\\
c_2(\zeta_{\eta},\eta) , & \text{if } y = 0,
\end{cases}
\end{equation}
where $J^*_{\eta-y}$ is the cost-to-go function that gives the remaining cost over the final $y$ time steps when applying an optimal series of inputs.  The principle of optimality \cite{bell66} implies that $J^*_{\eta-y}$ can be found iteratively according to
\begin{equation} \label{iterativecost}
J^*_{\eta-y} (\zeta_{\eta-y}) = \min_{u_{\eta-y}}  \bigg(  c_1(u_{\eta-y},\eta-y) + c_2(\zeta_{\eta-y},\eta-y)  +  J^*_{\eta-y+1}( \zeta_{\eta-y+1} ) \bigg),
\end{equation}
by starting with the endpoint cost $J^*_\eta(\zeta_\eta) = c_2(\zeta_{\eta}, \eta)$ and working backwards.  With knowledge of the cost-to-go function, the optimal control at each timestep is then
\begin{equation}  \label{optimalinput}
u_k^*(\zeta_k) = \underset{u_k \subset \mathcal{U}}{\rm argmin} \bigg(   c_1( u _{k},k)   +  c_2(\zeta_k,k) + J^*_{k+1}( f(\zeta_k,u_k) )    \bigg).
\end{equation}

\end{appendices}


\end{document}